\documentclass[10pt]{article}
\usepackage{latexsym,amssymb,amsmath,bbm,amsfonts,amsthm,
mathrsfs,bbm,euscript,pifont,lipsum,mdframed}

\def\cvd{\hfill$\Box$}
\def\vuoto{\varnothing}

\title{\Large{On the Fitting length of finite soluble groups I} \\ \vspace{1mm} 
\mbox{ \large Hall subgroups}}

\author{Giorgio Busetto -- Enrico Jabara}

\date{}

\begin{document}

\maketitle{}

\noindent {\bf Abstract.}  {\it Let $G$ be a finite soluble group 
and $h(G)$ its Fitting length. The aim of this paper is to 
give certain upper bounds for $h(G)$ as functions of the Fitting 
length of at least three Hall subgroups of $G$  which factorize 
$G$ in a  particular way.}

\vspace{1mm}


\vspace{10mm}

\textbf{\S 1. Introduction.}

\vspace{3mm}

If $G$ is a group containing subgroups $A$ and $B$ such that 
\[ G=AB=\{ab \mid a \in A, \; b \in B\} \] we say that $G$ is factorized 
by $A$ and $B$.
A classical problem in group theory is, generally speaking, to try to 
get informations on the structure of $G$  as a consequence of the structure 
of $A$ and $B$. For example a celebrated theorem by It\^{o} (\cite{Ito}) 
states that any group $G=AB$ with $A$ and $B$ abelian is metabelian.

In this paper we only deal with finite soluble groups and so 
for us group will always mean ``finite soluble group''. 
We shall consider the case when a group is factorized by Hall 
subgroups, namely subgroups whose order and index are coprime.
    
If $G$ is a group, we shall denote by $d(G)$ and $h(G)$ 
respectively the derived length and the Fitting length of $G$, 
by  $\pi(G)$ the set of prime divisors of $\vert G \vert$ 
and by $w(G)$ the cardinality of $\pi(G)$; also, we shall write 
$\pi$ and $w$ instead of $\pi(G)$ and $w(G)$ respectively when 
there is no possible ambiguity. 

A recent result in \cite{CJS} states that if  
$G=AB$, where $A$ and $B$ are Hall subgroups of $G$, then
\[ h(G) \leq h(A)+ h(B)+ 4d(B)-1.\] 
Under the previous hypotheses it is not possible to give an 
upper bound for $h(G)$ as a function of $h(A)$ and $h(B)$ only; 
for, it is well known that there exist groups of arbitrary 
Fitting length which can be factorized by two Sylow subgroups.
The aim of this paper is to prove that if $A$, $B$, 
$C$ are Hall subgroups of a group $G$ such that $G$ 
is trifactorized by $A$, $B$, $C$, namely if
 \[G=AB=BC=CA,\]
then there exists an upper bound for $h(G)$ as a function of  
$h(A)$, $h(B)$ and $h(C)$.
Moreover, if $G$ is $k$-factorized by Hall subgroups ($k\geq3$; 
the definition of $k$-factorized is straightforward), 
such upper bound gets better as $k$ increases.

If $G$ is a group and $\sigma$ is a set of primes we shall 
denote by $G_{\sigma}$ a Hall subgroup of $G$ such that
$\pi(G_{\sigma}) \subseteq \sigma$; if 
$\sigma=\pi(G) \setminus \{p\}$ for a prime $p$, we shall write 
$G_{p'}$ instead of $G_{\sigma}$.

\vspace{3mm}

\textsc{Theorem A.} {\it Let $G$ be a group and let $\sigma,\tau, \upsilon$
be three subsets of $\pi(G)$ such 
$\sigma \cup  \tau = \tau \cup \upsilon=\upsilon \cup 
\sigma=\pi$.  Then
\[ h(G) \leq h(G_{\sigma})+h(G_{\tau})+h(G_{\upsilon})-2\]
In particular, if $p,q \in \pi(G)$ and $p \not=q$, then}
\[ h(G) \leq h(G_{p'})+h(G_{q'})+h(G_{\{p,q\}})-2.\]

\vspace{3mm}

\textsc{Theorem B.} {\it Let $G$ be a group such that $w \geq 4$. 
Then there exist $p,q \in \pi(G)$  with $p \not =q$ such that}
\[ h(G) \leq h(G_{p'})+h(G_{q'})-1.\]

\vspace{3mm}

A simple example will show that Theorem B does not hold when $w=3$.

If  $G$  is a group we  define
\[ \mathfrak{h}_{\ell}(G)=\max \big\{ h(G_{\sigma}) \mid \sigma 
\subseteq \pi(G),\; \vert \sigma \vert=\ell \, \big\}\]
and we shall write $\mathfrak{h}_{\ell}$  instead of 
$\mathfrak{h}_{\ell}(G)$ if there is no possibile ambiguity.

\vspace{3mm}

\textsc{Theorem C.}  {\it Let $G$ be a group such that  $w \geq 3$ and 
assume $3 \leq \ell \leq w$. Then  
\[ \mathfrak{h}_{\ell} \leq 
\frac{\ell \cdot \mathfrak{h}_{\ell-1}-2}{\ell-2}.\]}

\vspace{3mm}

From Theorem C  we deduce
\[ h(G) \leq \frac{w \cdot \mathfrak{h}_{w-1}-2}{w-2}.\]
In particular, if $w =3$, then  
$h(G) \leq 3 \cdot \mathfrak{h}_{2}-2$ 
while if $w \geq 4$ then we have  $h(G) \leq 2 \cdot 
\mathfrak{h}_{w-1}-1$. The last bound is comparable with 
the one obtained in Theorem B.

\vspace{5mm}

\textbf{\S 2. Notation and preliminary results.}

\vspace{3mm}

As already pointed out, $G$ will always denote a finite 
soluble group and  $\{G_{p}\}_{p\in \pi}$ a Sylow system for $G$, 
namely a set of Sylow subgroups of $G$, one for any $p \in \pi$, 
such that $G_{p}G_{q}=G_{q}G_{p}$ for every $p,q \in \pi$.

If $\sigma \subseteq \pi$ by $\sigma$-Hall subgroup of $G$, 
we mean $G_{\sigma}=\prod_{p\in \sigma} G_{p}$; 
by $G_{p'}$ we denote the $\pi \setminus \{p\}$-Hall subgroup of $G$.
If $p \not \in \pi$, by definition  we set,
$G_{p'}=G$; if $\pi \subseteq \sigma$ then $G_{\sigma}=G$ and
if $\pi \cap \sigma=\vuoto$ then $G_{\sigma}=1$.

The rest of the notation will be mostly standard. 
In particular, if $p$ is a prime, $O_{p}(G)$ is the largest 
normal $p$-subgroup of $G$. We can write the Fitting subgroup 
$F(G)$ of $G$ as $F(G)= \prod_{p \in \pi} O_{p}(G)$.

We shall often make use of the following result (Lemma I.7.3 of
\cite{FT}; see also Lemma 1.2 of \cite{SS}).

\vspace{3mm}

\textsc{Lemma 2.1.} (\cite{FT}) {\it Let $G$ be a group and 
let $p \in \sigma \subseteq \pi$. 
Then for any $q \in \sigma$ such that
$q \not =p$ we have} \[ O_{p}(G_{\sigma}) \leq O_{q'}(G).\]

\vspace{3mm}

Using Lemma 2.1 we can deduce Lemma 1.3 of \cite{SS}, 
which we will prove in a way that allows us to get a better result.

\vspace{3mm}

\textsc{Proposition 2.2.}  {\it Let $G$  be a group and 
$p,q \in \pi$ with $p \not =q$. If $h(G_{p'}) \leq s$, 
$h(G_{q'}) \leq s$ and $h(G_{\{p,q\}}) \leq r$, then 
$h(G) \leq s(r+1)$.}

\vspace{1mm}

\textsc{Proof.} We use induction on $r$. If $r=0$ the statement 
is obvious. So, assume $r \geq 1$.
By Lemma 2.1 we have
\[ F(G_{\{p,q\}})=O_{p}(G_{\{p,q\}}) \times O_{q}(G_{\{p,q\}}) 
\leq O_{p'}(G)O_{q'}(G).\] 
Let $\overline{G}=G/\big(O_{p'}(G)O_{q'}(G)\big)$, then we have 
$h(\overline{G}_{\{p,q\}}) \leq r-1$ and, considering that  
$h\big(O_{p'}(G)O_{q'}(G)\big) \leq s$, the conclusion follows. \cvd

\vspace{3mm}

We shall not make use of Proposition 2.2 in what follows. 
It is presented here only to show how our Theorem A improves 
every upper bound  of this kind previously known in the literature.

In order to state our results we need the following 

\vspace{3mm}

\textsc{Definition 2.3.}  Let $G$ be a group, 
$t \geq 3$ be an integer and let $\mathscr{R}$ be the set
\[\mathscr{R}=\big\{ \varrho_{1}, \varrho_{2}, 
\ldots, \varrho_{t}  \mid \varrho_{i} \subseteq \pi 
\big \}. \]

Then $\mathscr{R}$ is called a \textit{cover} 
of order $t$ of $\pi$, shortly a $t$-cover, if
$\varrho_{i} \cup \varrho_{j} = \pi$
for every $i, j \in \{1,2, \ldots, t\}$, $i \not =j$.
A $t$-cover  is called  \textit{degenerate} 
if $\pi \in \mathscr{R}$.

The \textit{weight} of a $t$-cover $\mathscr{R}=
\{ \varrho_{1}, \varrho_{2},  \ldots, \varrho_{t}\}$ 
is the number
\[ \Theta(\mathscr{R})=\sum_{i=1}^{t} h(G_{\varrho_{i}}).\]

\vspace{3mm}

When no ambiguity is possible, we shall speak simply 
about cover  and we shall write
$\Theta$ for $\Theta(\mathscr{R})$.

\vspace{3mm}

\textsc{Remark 2.4} Let $G$ be a group, 
$\mathscr{R}$ a $t$-cover of $\pi$ and $w=w(G)$. 
Then
\begin{itemize}
\item[(a)] For every $p \in \pi$ there exists at most one 
$i \in \{1,2,\ldots, t\}$ such  that  $p \not \in \varrho_{i}$.
\item[(b)] The following inequality is always satisfied
\begin{equation}  \sum_{i=1}^{t} \vert \varrho_{i} \vert 
\geq (t-1) \cdot w. \end{equation}
\item[(c)] Suppose that $\mathscr{R}$
is a non degenerate $t$-cover of $\pi$. Then, by  (1)  it
follows  that  $t=\vert \mathscr{R} \vert \leq w$. 
Moreover, if $\vert \mathscr{R} \vert = w$, we have  $\vert
\varrho_{i} \vert = w -1$ for every $i \in \{1,2, \ldots, t\}$, 
namely for every $\varrho_{i}$ there exists $p_{i} \in \pi$ such
that $G_{\varrho_{i}}=G_{p_{i}'}$.
\item[(d)] If $\pi=\{p,q\}$, then $\pi$ admits a unique $3$-cover 
which is necessarily degenerate, namely 
$\mathscr{R}=\big\{ \{p\}, \{q\}, \pi \big\}$.
\end{itemize}

\vspace{3mm}

\textsc{Remark 2.5} Let $G$ be a group such that $w \geq 3$ 
and let $\mathscr{R}$ be a $t$-cover of $\pi$. 
Let $N\unlhd G$ and $\varrho_{i} \in \mathscr{R}$; if 
$\overline{G}=G/N$, then we can have
$\overline{G}_{\varrho_{i}}=\overline{G}$. 
In the proofs by induction on the order of $G$  sometimes we shall 
omit to mention explicitly this case. The reason for this is that we 
are interested in establishing an upper bound for $h(G)$ as a function 
of the values of $h(G_{\varrho_{i}})$. 
When $\overline{G}=\overline{G}_{\varrho_{i}}$ we have
$h(\overline{G})=h(\overline{G}_{\varrho_{i}})$ and  
the bound is automatically satisfied.

Similarly, if $\overline{G}_{\varrho_{i}}=1$, then 
there exists an index $j$ such that $\overline{G}_{\varrho_{j}}=\overline{G}$
and we are reduced to the case considered before. 

\vspace{5mm}

\textbf{\S 3.  Proofs, examples and final remarks.}

\vspace{3mm}

The following proposition is essential in order to get our 
results. The statement holds also when $w=2$ by Remark 2.4 (d); 
however, for our purposes, it is convenient to assume the hypothesis 
that $w \geq 3$.

\vspace{3mm}

\textsc{Proposition 3.1} {\it Let $G$ be a group such that $w\geq 3$  
and let $\mathscr{R}$ be $t$-cover
of $\pi$ of  weight $\Theta$. If $t \geq 3$, then} 
\begin{equation}  \boxed{ h(G) \leq \frac{\Theta-2}{t-2}}
 \end{equation}

\vspace{1mm}

\textsc{Proof.}  Let  $\mathscr{R}=\{\varrho_{1},\varrho_{2},
\ldots,\varrho_{t}\}$, $h_{i}=h(G_{\varrho_{i}})$. We prove the claim by 
induction on $t+\Theta+h(G)+\vert G \vert$. 

If $G$ is nilpotent, then  $\Theta=t$ and the result follows immediately. 
Hence we  may assume that $h(G) \geq 2$. 
Let $F=F(G)$ be the Fitting subgroup of $G$. We consider two cases:

\vspace{1mm} (\textbf{A}) $w(F) \geq 2$. \vspace{1mm} 

Then there exist $p,q \in \pi$ such that
$O_{p}(G) \not = 1 \not =O_{q}(G)$.
If neither $O_{p}(G)$ nor $O_{q}(G)$ 
are Sylow subgroups of $G$ , or if  $w(G) \geq 4$, then 
we set  $\overline{G}=G/O_{p}(G)$ and $\widehat{G}=G/O_{q}(G)$.
In $\overline{G}$ and $\widehat{G}$  
the hypotheses are preserved and a standard argument 
gives the conclusion.

Hence suppose that  $O_{p}(G)$ is a Sylow $p$-subgroup of  $G$ 
(clearly the roles of $p$ and $q$ are interchangeable) and that  
$w(G) =3$. 
Then $\pi(G)=\{p,q,r\}$ where $p$, $q$ and $r$ 
are different primes, and $t-2=1$.
If $\pi \in \mathscr{R}$, 
then the statement follows immediately. 
Therefore we may assume that $\mathscr{R}$ is not degenerate and 
$\mathscr{R}=\big\{ \{p,q\}, \{q,r\}, \{r,p\} \big\}$.
Note that $G/O_{p}(G) \simeq G_{q}G_{r}$.
If $h(G_{p}G_{q})=1=h(G_{p}G_{r})$, then
$G=G_{p} \times G_{q}G_{r}$ and
$h(G)=h(G_{q}G_{r})=h(G_{p}G_{q})+h(G_{q}G_{r})+
h(G_{r}G_{p})-2$. If $h(G_{p}G_{q})+
h(G_{p}G_{r}) \geq 3$, then
\[h(G) \leq h(G_{q}G_{r})+1 \leq 
h(G_{p}G_{q})+h(G_{q}G_{r})+h(G_{r}G_{p})-2.\]

\vspace{1mm} (\textbf{B}) $w(F)=1$ and $F=O_{p}(G)$
for some $p \in \pi$. \vspace{1mm} 

If $F$ is a Sylow $p$-subgroup of $G$ we can get the 
conclusion using Remark 2.5 and arguing as in (\textbf{A}). 
Hence we may assume that  $\pi(G/F) = \pi$.

By  Remark  2.4.(a), possibly reordering the indices,  
we may assume that  
$F \leq G_{\varrho_{i}}$ for $i=1,2, \ldots, t-1$.
Moreover, since $C_{G}(F) \leq F$ and
$C_{G_{\varrho_{i}}}(F(G_{\varrho_{i}})) 
\leq F(G_{\varrho_{i}})$, we deduce that $F(G_{\varrho_{i}})$ 
is a $p$-group containing $F$, namely 
$F(G_{\varrho_{i}})=O_{p}(G_{\varrho_{i}})$ for every
$i=1,2, \ldots, t-1$.

From Lemma 2.1 it follows that  
$O_{p}(G_{\varrho_{i}}) \leq O_{q_{i}'}(G)$ 
for every $q_{i} \in \varrho_{i} \setminus \{ p\}$.

Let us fix $q_{1} \not \in \varrho_{1}$;  
by  Remark 2.4.(a) we get $q_{1} \in \varrho_{i}$ 
per $i=2,3, \ldots, t-1$. 
In $\overline{G}=G/ O_{q_{1}'}(G)$ we have
$h(\overline{G}_{\varrho_{i}}) \leq h_{i}-1$ 
for every $i=2,3, \ldots t-1$ and, by inductive hypothesis 
\[ h(\overline{G}) \leq \frac{\left(\sum_{i=1}^{t} 
h(\overline{G}_{\varrho_{i}})\right) -2}{t-2} \leq 
\frac{\big(\Theta -t+2 \big)-2}{t-2}
= \frac{\Theta-2}{t-2} -1.\]

Similarly, if $q_{2} \not \in \varrho_{2}$, then
$q_{2} \in \varrho_{i}$ for $i=1,3, \ldots, t-1$. 
In $\widetilde{G}=G/ O_{q_{2}'}(G)$ we have
$h(\widetilde{G}_{i}) \leq h_{i}-1$ for
$i=1,3, \ldots t-1$ and, arguing as before, we obtain
\[h(\widetilde{G}) \leq \frac{\Theta-2}{t-2} -1.\]

\vspace{1mm}

By hypothesis $t-1 \geq 2$, hence
\[ \left(\bigcap_{r \not \in \{p,q_{1}\}} O_{r'}(G) \right) 
\cap \left( \bigcap_{r \not \in \{p,q_{2}\}} O_{r'}(G) \right) 
= \bigcap _{r \not =p} O_{r'}(G) = O_{p}(G), \]
therefore
\[ h(G)-1=h\big(G/O_{p}(G) \big) \leq \frac{\Theta-2}{t-2} -1,\]
and the conclusion follows.   \cvd  

\vspace{3mm}

If $t=3$ from Proposition 3.1 it follows that $h(G) \leq \Theta-2$;
if $G$ is nilpotent such a bound is obviously  the best. 
We give some examples showing that the bound given by (2) is 
sufficiently accurate for any value of $h(G)$.

In the Examples 3.2, 3.3, 3.4 and 3.5  
we  shall make use of the following notation.
If $H$ is a group and $\ell \geq 1$, we shall write
$[H]_{\ell}$ for the  iterated wreath product of $H$, namely 
$[H]_{1}=H$ and $[H]_{\ell+1}=[H]_{\ell} \wr H$. 
If $p$, $q$ and $r$ are different prime numbers, 
let us denote by $P$, $Q$, $R$ a $p$-group, 
a $q$-group and an $r$-group respectively. 
If $\pi=\{p,q,r\}$, then we set
$\mathscr{R}=\big\{ \{p,q\}, \{q,r\},\{r,p\} \big\}$. 

\vspace{3mm}

\textsc{Example 3.2.}  Let $G=P\wr [Q \wr R]_{\ell}$. 
Then $h(G)=2\ell+1$, $h(G_{r'})=h(G_{q'})=2$ and
$h(G_{p'})=2\ell$, therefore $\Theta-2=2\ell+2$.

But, if $G=P \times [Q \wr R]_{\ell}$ then $h(G)=2\ell$,
$h(G_{r'})=h(G_{q'})=1$, $h(G_{p'})=2\ell$ and therefore 
$\Theta -2=2\ell$.

\vspace{3mm}

\textsc{Example  3.3.} Let $G=([P \wr Q]_{\ell}) \wr 
([R \wr Q]_{\ell})$. 
It can be seen that $h(G)=4\ell$ whereas $h(G_{r'})=2\ell$,
$h(G_{p'})=2\ell+1$ and $h(G_{q'})=2$. Hence $\Theta-2=4\ell+1$.

\vspace{3mm}

\textsc{Proof of Theorem  A.} With the notation of Proposition 3.1, 
if we set $t=3$, 
$\varrho_{1}=\pi \setminus \{p\}$, $\varrho_{2}=\pi \setminus 
\{q\}$ and $\varrho_{3}= \{p,q\}$, then the  conclusion follows. \cvd

\vspace{3mm}

For notational convenience, given the group $G$, 
in the following proofs we shall denote by  
$\mathscr{R}^{\ast}$ the $w$-cover 
$\big\{ \pi \setminus \{p\} \mid p \in \pi \big\}$ of $\pi(G)$.

\vspace{3mm}

\textsc{Proof of Theorem B.} Let $G$ be a group such that $w \geq 4$. 
We have to show that there exist $p,q \in \pi$ such that 
$h(G) \leq h(G_{p'})+h(G_{q'})-1$. 

Let $\pi=\{ p_{1},p_{2}, \ldots, p_{w} \}$, 
where the indices are ordered in a way that  \[ h(G_{p_{1}'}) 
\geq  h(G_{p_{2}'}) \geq \ldots \geq h(G_{p_{w}'}). \]
Set $p=p_{1}$ and $q=p_{2}$. 

If $h(G_{p'})=1$ then $G$ is nilpotent and the statement holds. 
If $h(G_{p'}) \geq 2$ then we cannot have $h(G_{q'})=1$;  
for, $w \geq 4$ and from Theorem A it would follows that
$h(G) \leq h(G_{p_{2}'})+h(G_{p_{3}'})+h(G_{p_{4}'})
-2=1$, hence $h(G)=1$.

If we set $\lambda= h(G_{p'})+h(G_{q'})$ then we may assume that  
$\lambda \geq 4$.
Considering the hypothesis that $w \geq 4$, by easy 
considerations we get
\begin{equation} 
w \lambda -4 \leq 2  (w-2)  (\lambda-1). 
\end{equation}

Given the cover $\mathscr{R}^{\ast}$, we certainly have
\[ \Theta^{\ast} \leq \frac{w}{2} \cdot 
\big(h(G_{p'})+h(G_{q'}) \big) =\frac{w  \lambda}{2}. \]
Therefore, by applying the inequality (3) 
and Proposition 3.1, we get
\[ h(G) \leq \frac{\Theta^{\ast}-2}{w-2} 
\leq \frac{w\lambda -4}{2(w-2)} \leq \lambda -1=h(G_{p'})+h(G_{q'}) -1,\]
and the statement is proved. \cvd

\vspace{3mm}

It is not possible to extend Theorem B to the case $w=3$; 
for, there exist the following

\vspace{3mm}

\textsc{Example 3.4.} Let $G=[P \wr Q]_{\ell} \wr  [R \wr P]_{\ell} 
\wr [Q \wr R]_{\ell}$ where the order of the wreath products can be 
chosen as you prefer. We have   
$h(G_{p'})=h(G_{q'})=h(G_{r'})=2\ell+2$, while $h(G)=6\ell$.

\vspace{3mm}

\textsc{Example 3.5.} Let $s$ be a prime, $s \not \in \{p,q,r \}$ 
and let $S$ be a non-trivial $s$-group. 
Let\[ G=[P\wr Q]_{\ell} \wr [P \wr R ]_{\ell} \wr [P \wr S]_{\ell} 
\wr [Q \wr S]_{\ell} \wr [Q \wr P]_{\ell} \wr [Q \wr R]_{\ell}\]
where the order of the wreath products can be chosen as you prefer.
We have  $h(G)=12\ell$, $h(G_{p'})=4\ell+3$, $h(G_{q'})=4\ell+2$
and $h(G_{\{p,q\}})=4\ell+2$ hence $h(G_{p'})+h(G_{q'})+h(G_{\{p,q\}})-2
=12\ell +5$. Finally \[ 
\frac{\Theta(\mathscr{R}^{\ast})-2}{4-2}=\frac{1}{2}\big( 
h(G_{p'})+h(G_{q'})+h(G_{r'})+h(G_{s'})-2\big)=12\ell+2.\]

\vspace{3mm}

\textsc{Proof of Theorem C.} Let $G$ be a group such that
$w \geq 3$ and let $\ell \in \mathbb{N}$ where
$3 \leq \ell \leq w$.  We show that
\begin{equation} \mathfrak{h}_{\ell} \leq 
\frac{\ell \cdot \mathfrak{h}_{\ell-1}-2}{\ell-2}.
\end{equation}

Let $H$ be a Hall subgroup of $G$  such that 
$w(H)=\ell$ and $h(H)=\mathfrak{h}_{\ell}$. 
In order to prove the statement it is not restrictive to 
assume $G=H$. Given the cover $\mathscr{R}^{\ast}$, 
then, for every $\varrho \in \mathscr{R}^{\ast}$, we
have $h(G_{\varrho}) \leq \mathfrak{h}_{\ell-1}$ and 
therefore $\Theta^{\ast} \leq \ell \cdot \mathfrak{h}_{\ell-1}$. 
The conclusion follows by Proposition 3.1. \cvd

\vspace{3mm}

\textsc{Remark 3.6.} Let $G$ be a group such that $w \geq 3$.  
Using  induction  and applying the inequality (4), we obtain
\[ h(G)=\mathfrak{h}_{w} \leq  \frac{w(w-1)}{2} \cdot 
\big(\mathfrak{h}_{2}-1 \big) +1. \]

\vspace{3mm}

We observe that it seems quite difficult to obtain our 
results omitting the hypothesis that the subgroups 
considered are Hall subgroups. 
In this direction we make the following 

\vspace{3mm}

\textsc{Conjecture 3.7.}  Let $G$ be a group and suppose 
that $H, K, L \leq G$  are subgroups of  $G$  such that
$G=HK=KL=LH$. 
Then
\begin{equation} h(G) \leq h(H)+h(K)+h(L)-2 \end{equation} 
or, at least, there exists a function $\eta$ such that 
$h(G) \leq \eta \big(h(H),h(K),h(L)\big)$.

\vspace{3mm}

The unique case we know where the inequality (5)  holds is when 
$h(H)=h(K)=1$. For, if 
$H$, $K$ and $L$ are nilpotent, then a result by Kegel 
(\cite{Ke}, see also Corollary 2.5.11 of \cite{AFG}), states that $G$ 
is also nilpotent. The general case is covered by an unpublished result 
by Peterson (see Theorem 2.5.10 of \cite{AFG}).
Probably, in order to obtain some result in this new direction, 
it is necessary to use hypotheses stronger than the ones considered 
in Conjecture 3.7. Hence we also make the following

\vspace{3mm}

\textsc{Conjecture 3.8.} Let $G$ be a group and suppose that 
$N_{1}$, $N_{2}$, $N_{3}$ are nilpotent subgroups of $G$
such that $G=N_{1}N_{2}N_{3}$ and $N_{i}N_{j}=N_{j}N_{i}$
for every $i,j \in \{1,2,3\}$. Then
\[ h(G) \leq h(N_{1}N_{2})+h(N_{2}N_{3})+h(N_{3}N_{1})-2.\]

\vspace{8mm}

\textsc{Giorgio Busetto}

\textsc{DAIS -- Universit\`{a} di Venezia}

\textsc{Dorsoduro 2137, 30123 Venezia - ITALY}

\texttt{gbusetto@unive.it}

\vspace{3mm}

\textsc{Enrico Jabara}

\textsc{DFBBC -- Universit\`{a} di Venezia}

\textsc{Dorsoduro 3484/D, 30123 Venezia - ITALY}

\texttt{jabara@unive.it}

\end{document}